\documentstyle[12pt,amssymb]{article}
\begin{document}
\def\y {\'{\i}}
\def\beq{\begin{equation}}
\def\eeq{\end{equation}}
\def\bray{\begin{eqnarray}}
\def\eray{\end{eqnarray}}

\begin{titlepage}
${}$
March 1996 \hfill{IFT-P.007/96}

\hfill{q-alg/9603004}
\vskip .6in
 
\begin{center}
{\large {\bf Non Abelian Sugawara Construction and the q-Deformed N=2
Superconformal Algebra}}
\end{center}

\normalsize
\vskip .4in
 
\begin{center}

{E. Batista \footnotemark \footnotetext{Supported by FAPESP} ,
J. F. Gomes\footnotemark \footnotetext{Work partially
 supported by CNPq} , 
I.J. Lautenschleguer \footnotemark \footnotetext {Supported by CNPq}} \\
\par \vskip .1in \noindent
Instituto de F\'{\i}sica Te\'{o}rica-UNESP\\
Rua Pamplona 145\\
01405-900 S\~{a}o Paulo, S.P.\\
 Brazil \\
\par \vskip .3in

\end{center}
 
\begin{center}
{\large {\bf ABSTRACT}}\\
The construction of a q-deformed N=2 superconformal algebra is proposed in
terms of level 1 currents of ${\cal{U}}_{q} ({\widehat{su}}(2))$ quantum
 affine 
Lie algebra and a single real Fermi field. In particular, it suggests the
expression for the q-deformed Energy-Momentum tensor in the Sugawara form.
  Its 
constituents  generate two isomorphic quadratic algebraic structures. The
generalization to ${\cal{U}}_{q} ({\widehat{su}}(N+1))$ is also
proposed.

\end{center}
\par \vskip .3in \noindent

\end{titlepage}

\section{Introduction}

It has become quite well known that conformal invariance uncovers a common
structure among many important field theoretical models in two dimensions.  
The physics content of such models may be extracted by use of representation
theory of the underlying conformal algebra.  String theory is among the many
important examples where the Energy-Momentum tensor is written in the
Sugawara form, i.e. bilinear in conserved currents.

Quantum groups has also played an important role in the sense that a class of 
perturbed systems may be casted within a deformed algebraic structure
providing, by representation theory, the physics content of the model. This
algebraic structure is characterized by a deformation parameter q.  In
particular, a q-deformed version of the Veneziano model was proposed in
\cite{cgk} and \cite{cgr} by replacing the ordinary oscillators by q-deformed
ones in the operator formalism of Fubini and Veneziano.  The deformed model
(when $q\neq 1$) was shown to lead to non linear Regge trajectories. It is hoped
that a conformal structure arises in terms of a q-deformed Kac-Moody algebra
such that the usual Sugawara construction is smoothly recovered in the limit
$q\rightarrow 1$.  Also in this line, a q-deformed proposal for the Nambu
action was discussed by de Vega and Sanchez in \cite{devega}.  Their action 
shows a nonlocal character and the canonical Energy-Momentum tensor
 was proposed bilinearly in terms of q-oscillators.

In this paper we study how a q-conformal structure may be constructed in terms
of the oscillators proposed by Frenkel and Jing \cite{fjing}. They have
constructed vertex operators for simply laced Lie algebras satisfying a level 1
q-deformed Kac-Moody algebra under the operator product expansion (OPE). 

An interesting feature of a q-deformed field theory is that the OPE's turns out
to be less divergent than the $q = 1$ theory.  This fact shows up since poles
are smeared out in a symmetric manner in terms of the deformation parameter
$q$  signaling towards a nonlocal structure.  The
short distance behaviour of the operators require a  consistent definition of
normal ordering. In order to produce the correct analytic structure in terms
of poles we shall consider the deformation parameter restricted to a pure phase.

In section 2 we discuss and review the abelian q-Sugawara construction and the
$N=1$ q-superconformal algebra obtained from a 
level 1  q-deformed $H_q (\infty )$ infinite Heisenberg algebra and a
single real fermi field.  In section 3 we discuss the vertex operator
construction of Frenkel and Jing for ${\cal{U}}_q ({\widehat{su}} (2))$ 
and derive in appendix many useful identities using the q-Taylor expansion.

In section 4, guided by the basic principle of decomposing the most divergent
poles into a product of simple poles shifted in a symmetric manner in terms of
the deformation parameter $q$,  together with the closure of algebra, we 
were led to define a family of super charge generators $G_{\alpha}^{\pm}$ 
in terms of the vertex operators of ${\cal{U}}_q ({\widehat{su}} (2))$ and a
single real fermi field.  We show that this construction leads to a q-deformed
$N=2$ superconformal algebra.  In particular, the OPE of the two super charge
generators define the Energy-Momentum tensor which is shown to decompose
into commuting  bosonic and fermionic couterparts as forcasted by 
Goddard and Schwimmer
\cite{goddard}. The fermionic counterpart of the Energy-Momentum
tensor defines a q-deformed Virasoro algebra as in ref. \cite{sato1}. 
The bosonic Energy-Momentum tensor however is no longer bilinear in
q-oscillators and does not yield a closed algebra 
but it can be written as a sum of two components, each of which generating
 closed algebras of quadratic type. Similar structures have
recently been constructed by Frenkel and Reshetikhin \cite{fr} 
using Wakimoto
realization for ${\cal{U}}_q ({\widehat{sl}} (2))$ currents. Finally, we
generalize, in the last section, the results for other simply laced q-deformed 
affine Lie algebras.

\section{The Abelian Sugawara Construction and The 
q-Deformed $N=1$ Superconformal Algebra}

An abelian q-deformed Kac-Moody algebra can be constructed from the usual
undeformed level 1 $U (1)$ current algebra,
\beq
H(z) H(w) \, =\,  \frac{1}{(z\, -\, w)^{2}}  +
  Regular \; Terms \;   
\label{u1}
\eeq
by replacing the double pole into a product of two simple poles
 symmetrically displaced in terms of a deformation parameter $q$, i.e.
\beq
\label{umope}
H(z) H(w) \, =\, \frac{1}{(z - qw)(z-wq^{-1})}  
+  Regular \; Terms.
\eeq

The abelian field $H(z)$ can be Laurent expanded in terms of a
q-deformed infinite Heisenberg algebra $H_q (\infty )$
\begin{equation}
\label{umgenfunction}
H(z) \, =\, \sum_{n\in \Bbb Z} a_n z^{-n-1} \quad .
\end{equation}
where from (\ref{umope}) we find 
\begin{equation}
\label{umosc}
[a_n ,a_m] \, =\, [n] \, \delta_{m+n,0} \quad ,
\end{equation}
where  $[n] \, =\, \frac{q^n -q^{-n}}{q-q^{-1}}$.
Because the non-local character of the pole structure, to obtain a consistent
definition of normal ordering, we shall be assuming the deformation parameter 
to be a pure phase, i.e. $q=e^{i\epsilon }\, ,\, \epsilon \in \Bbb R$.
The Sugawara construction implies the energy momentum tensor to be written as
\begin{equation}
\label{suga}
T(z) \, =\, \frac{1}{2} :H(z) H(z) : \quad ,
\end{equation}
where the dots denote normal ordering in the sense that the positive oscillator
modes are moved to the right of those with negative. A classical version of
(\ref{suga}) was proposed in ref \cite{chai} while normal ordering was later 
introduced in \cite{singh} with the effect of generating central terms in the
q-Virasoro algebra.  From (\ref{umope}), it is straightforward to evaluate
\beq
\label{ht}
T(z) H(w) = 
\frac{1}{w(q - q^{-1})} \left (
\frac{H(wq)}
{z - wq}  - 
\frac{H(wq^{ -1})}
{z - wq^{-1}} \right ) \, + Reg \;  Terms, \;  |z|>|w|
\eeq

The shift in the arguments of $H(w)$ in the r.h.s. requires an extra index for
the closure of the algebra of the energy momentum tensor (\ref{suga}).   Define
\begin{equation}
\label{umemboson}
T^{\alpha}(z) \, =\, \frac{1}{2} :H(zq^{\frac{\alpha}{2}}) 
H(zq^{-\frac{\alpha}{2}}) : \quad ,\, \alpha \in \Bbb Z \quad ,
\end{equation}
It them follows from (\ref{ht}) that these operators satisfy a deformation of
the Virasoro algebra with unit central charge, 
\bray
T^{\alpha}(z) T^{\beta}(w) \, &=&\, 
\frac{1}{wq^{\frac{\beta +\alpha}{2}}(q\, -\, q^{-1})} \left (
\frac{T^{\beta -\alpha +1}(wq^{\frac{-\alpha +1}{2}})}
{z\, -\, wq^{\frac{\beta -\alpha}{2}+1}} \, -\, 
\frac{T^{\beta -\alpha -1}(wq^{\frac{-\alpha -1}{2}})}
{z\, -\, wq^{\frac{\beta -\alpha}{2}-1}} \right ) \,  
\nonumber \\
&+&\,
\frac{1}{wq^{\frac{-\beta +\alpha}{2}}(q\, -\, q^{-1})} \left (
\frac{T^{-\beta -\alpha +1}(wq^{\frac{-\alpha +1}{2}})}
{z\, -\, wq^{\frac{-\beta -\alpha}{2}+1}} \, -\, 
\frac{T^{-\beta -\alpha -1}(wq^{\frac{-\alpha -1}{2}})}
{z\, -\, wq^{\frac{-\beta -\alpha}{2}-1}} \right ) \,  
\nonumber \\
&+&\, 
\frac{1}{wq^{\frac{\beta -\alpha}{2}}(q\, -\, q^{-1})} \left (
\frac{T^{\beta +\alpha +1}(wq^{\frac{\alpha +1}{2}})}
{z\, -\, wq^{\frac{\beta +\alpha}{2}+1}} \, -\, 
\frac{T^{\beta +\alpha -1}(wq^{\frac{\alpha -1}{2}})}
{z\, -\, wq^{\frac{\beta +\alpha}{2}-1}} \right ) \,  
\nonumber \\
&+&\, 
\frac{1}{wq^{\frac{-\beta -\alpha}{2}}(q\, -\, q^{-1})} \left (
\frac{T^{-\beta +\alpha +1}(wq^{\frac{\alpha +1}{2}})}
{z\, -\, wq^{\frac{-\beta +\alpha}{2}+1}} \, -\, 
\frac{T^{-\beta +\alpha -1}(wq^{\frac{\alpha -1}{2}})}
{z\, -\, wq^{\frac{-\beta +\alpha}{2}-1}} \right ) \,  
\nonumber
\eray
\[
+\, \frac{1}{4(z\, -\, wq^{\frac{\beta -\alpha}{2}+1})
(z\, -\, wq^{\frac{\beta -\alpha}{2}-1})  
(z\, -\, wq^{\frac{-\beta +\alpha}{2}+1})
(z\, -\, wq^{\frac{-\beta +\alpha}{2}-1})} \,  
\]
\begin{equation}
\label{umviraboson}
+\, \frac{1}{4(z\, -\, wq^{\frac{-\beta -\alpha}{2}+1})
(z\, -\, wq^{\frac{-\beta -\alpha}{2}-1}) 
(z\, -\,  wq^{\frac{\beta +\alpha}{2}+1})
(z\, -\,  wq^{\frac{\beta +\alpha}{2}-1})} \quad   |z|>|w|
\end{equation}

A construction of a supercharge generator was proposed in \cite{chai} by
introducing a real fermi field $\psi (z)$ 
\begin{equation}
\label{umfermiope}
\psi (z) \psi (w) \, =\, :\psi (z) \psi (w) : \, +\, \frac{1}{(z\, -\, w)}
\quad  |z|>|w|.
\end{equation}
Since the triple pole produced by the product of two supercharge generators
in the $q=1$ case should be replaced by  the product of three simple 
poles, 
\[\frac{1}{(z-w)^3} \rightarrow \frac{1}{(z-wq^{-1})(z-w)(z-wq)} \quad ,
\]
we realize this by the product of the
bosonic current (\ref{umope}) and the fermion (\ref{umfermiope}), ie, 
$G(z)=H(z) \psi(z)$. The closure condition of the algebra requires 
  a family of supercharge generators to be defined as
\begin{equation}
\label{umsupercarga}
G^{\alpha} (z) \, =\, H(zq^{\frac{\alpha}{2}}) \psi(zq^{-\frac{\alpha}{2}})
\quad .
\end{equation}
Their algebra is given by
\bray
G^{\alpha}(z) G^{\beta}(w) \, &=&\,  
\frac{2q^{\frac{\alpha}{2}} T^{\beta +\alpha}(wq^{\frac{\alpha}{2}})}
{z\, -\, wq^{\frac{-\beta +\alpha}{2}}} \, \nonumber \\
&+&\, 
\frac{[\beta -\alpha +1] L^{\beta -\alpha +1} (wq^{\frac{-\alpha +1}{2}})}
{q^{\frac{\beta -1}{2} +\alpha}(z\, -\, wq^{\frac{\beta -\alpha}{2} +1})} \, 
\nonumber \\
&-&\, 
\frac{[\beta -\alpha -1] L^{\beta -\alpha -1} (wq^{\frac{-\alpha -1}{2}})}
{q^{\frac{\beta +1}{2} +\alpha}(z\, -\, wq^{\frac{\beta -\alpha}{2} -1})} 
\,  
\nonumber
\eray
\begin{equation}
\label{umsuperesuper}
+\, \frac{q^{-\frac{\alpha}{2}}}
{(z\, -\, wq^{\frac{-\beta +\alpha}{2}})
(z\, -\, wq^{\frac{\beta -\alpha}{2}+1})
(z\, -\, wq^{\frac{\beta -\alpha}{2}+1})}   \quad  |z|>|w|.
\end{equation}
The operators $L^{\alpha} (z)$ are the fermionic part of the Energy-Momentum
tensor
\begin{equation}
\label{umemfermi}
L^{\alpha} (z) \, =\, \frac{1}{[\alpha](q\, -\, q^{-1})}
 :\psi (zq^{\frac{\alpha}{2}}) \psi (zq^{-\frac{\alpha}{2}}): \quad 
\end{equation}
and satisfy a q-deformed Virasoro algebra (see \cite{sato1}) 
\bray
\label{umvirafermi}
L^{\alpha} (z) L^{\beta} (w) \, &=&\, 
\frac{1}{[\alpha ][\beta ]w(q\, -\, q^{-1})} \left (
\frac{[\alpha +\beta]L^{\alpha +\beta}(wq^{\frac{\alpha}{2}})}
{q^{\frac{\beta -\alpha}{2}}(z\, -\, wq^{\frac{\beta +\alpha}{2}})} \, 
\right. \nonumber \\
&+&\, 
\frac{[-\alpha -\beta]L^{-\alpha -\beta}(wq^{-\frac{\alpha}{2}})}
{q^{\frac{-\beta +\alpha}{2}}(z\, -\, wq^{\frac{-\beta -\alpha}{2}})} \, -\,
\frac{[\alpha -\beta]L^{\alpha -\beta}(wq^{\frac{\alpha}{2}})}
{q^{\frac{-\beta -\alpha}{2}}(z\, -\, wq^{\frac{-\beta +\alpha}{2}})} \, 
\nonumber \\
&-&\,
\frac{[-\alpha +\beta]L^{-\alpha +\beta}(wq^{\frac{-\alpha}{2}})}
{q^{\frac{\beta +\alpha}{2}}(z\, -\, wq^{\frac{\beta -\alpha}{2}})} \, \,
+ \frac{1}
{(z\, -\, wq^{\frac{-\beta -\alpha}{2}})(z\, -\, wq^{\frac{\beta +\alpha}{2}})}
\nonumber \\
&-&\, 
\left. \frac{1}
{(z\, -\, wq^{\frac{-\beta +\alpha}{2}})(z\, -\, wq^{\frac{\beta -\alpha}{2}})}
\right ) \quad  |z|>|w|.
\eray
This algebra is not isomorphic to the bosonic one given in (\ref{umviraboson}),
 but in the limit $q\rightarrow 1$, these generators become the usual fermionic
Energy-Momentum tensor
\[
L(z) \, =\, \frac{1}{2} :\partial_z \psi (z) \psi (z): \quad ,
\]
satisfying the usual Virasoro algebra with $c=\frac{1}{2}$.

The closure of the $N=1$ super conformal algebra is achieved with the following
OPE relations
\bray
\label{umbosonesuper}
T^{\alpha}(z) G^{\beta}(w)\, &=&\, 
\frac{1}{2wq^{\frac{\beta}{2}} (q\, -\, q^{-1})} \left (
\frac{q^{-\frac{\alpha}{2}} G^{\beta -\alpha +1} (wq^{\frac{-\alpha +1}{2}})}
{z\, -\, wq^{\frac{\beta -\alpha}{2} +1}} \, \right. \nonumber\\
&-&\, 
\frac{q^{-\frac{\alpha}{2}} G^{\beta -\alpha -1} (wq^{\frac{-\alpha -1}{2}})}
{z\, -\, wq^{\frac{\beta -\alpha}{2} -1}} \, + \, 
\frac{q^{\frac{\alpha}{2}} G^{\beta +\alpha +1} (wq^{\frac{\alpha +1}{2}})}
{z\, -\, wq^{\frac{\beta +\alpha}{2} +1}} \, \nonumber\\
&-&\, \left. 
\frac{q^{\frac{\alpha}{2}} G^{\beta +\alpha -1} (wq^{\frac{\alpha -1}{2}})}
{z\, -\, wq^{\frac{\beta +\alpha}{2} -1}} \right ) \quad  |z|>|w|,
\eray
\begin{equation}
\label{umfermiesuper}
L^{\alpha}(z) G^{\beta}(w) \, =\, 
\frac{1}{[\alpha]wq^{-\frac{\beta}{2}}(q\, -\, q^{-1})} \left (
\frac{G^{\beta -\alpha}(wq^{\frac{\alpha}{2}})}
{z\, -\, wq^{\frac{-\beta +\alpha}{2}}} \, -\,  
\frac{G^{\beta +\alpha}(wq^{\frac{-\alpha}{2}})}
{z\, -\, wq^{\frac{-\beta -\alpha}{2}}} \right ) \quad  |z|>|w|,
\end{equation}
In the limit $q \rightarrow 1 $ we recover the usual undeformed $N=1$
superconformal algebra, namely
\bray
G(z) G(w) \, &=&\, \frac{1}{(z\, -\, w)^3} \, + \frac{2 T(w)}{(z\, -\, w)}
\quad , \nonumber \\
T(z) G(w) \, &=&\, \frac{\frac{3}{2} G(w)}{(z\, -\, w)^2} \, +\, 
\frac{G'(w)}{(z\, -\, w)} \quad , \nonumber \\
T(z) T(w) \, &=&\, \frac{2T(w)}{(z\, -\, w)^2} \, +\, 
\frac{T'(z)}{(z\, -\, w)} \, +\, \frac{\frac{3}{4}}{(z\, -\, w)^4} \quad .
\nonumber 
\eray
where $T = T_{bosonic} + L_{fermionic}$.

The deformation is therefore responsible for introducing an addicional 
integer  index
into the algebraic structure.  Similar structures were already been found in
the context of two loop Kac-Moody algebras and with the conformal affine Toda 
models (see for instance refs. \cite{afgz},\cite{fgsz} ).
 We should also point out that the decomposition
  of higher order  poles into the product of simple
poles does not alter the merororphic character of the theory.  In fact, a
consequence of this is the splitting of the bosonic and fermionic parts of the
total energy momentum tensor ( see \cite{goddard}).

\section{ Deformed ${\cal{U}}_{q}({\widehat{su}}(2))$ Kac-Moody Algebras
 and Vertex Operators} 
                                           
Consider a level 1  ${\cal{U}}_{q}({\widehat{su}}(2))$ construction 
in terms of vertex operators proposed by Frenkel and Jing \cite{fjing} 
using the q-deformed oscillators 
\begin{equation}
\label{2osc}
[\alpha_n ,\alpha_m ] \, =\, \frac{[2n][n]}{2n} \, \delta_{m+n,0} \quad ,
\end{equation}
together with the undeformed Heisenberg algebra       
\begin{equation}            
\label{zeromodes}
[Q,P] \, =\, i  \quad .
\end{equation}
Define a Fubini-Veneziano field 
\bray
\label{2fubini}
\xi^{\pm} (z) \, &=&\, Q \, -\, iP\ln{z} \, +\, i\sum_{n<0}
\frac{\alpha_n}{[n]} (zq^{\mp \frac{1}{2}})^{-n} \,  \nonumber \\
&+&\, 
i\sum_{n>0} \frac{\alpha_n}{[n]} (zq^{\pm \frac{1}{2}})^{-n} \quad .
\eray
and then write the generators for  ${\cal{U}}_{q}({\widehat{su}}(2))$ as
\begin{equation}
\label{2vertex}
E^{\pm} (z) \, =\, :\exp \left\{ \pm \sqrt{2} \xi^{\pm} (z) \right\} : \quad ;
\, H(z) \, =\, \sum_{n\in \Bbb Z} \alpha_n z^{-n-1} \quad ,
\end{equation}
where the dots denote normal ordering in the sense that positive oscillator
modes are moved to the right of those with negative and P to the right of Q.
The OPE for vertices like (\ref{2vertex}) is obtained using the
 Baker-Campbell-Haussdorff formula.  In particular we find,
\begin{equation}
\label{2baker}
E^{+} (z) E^{-} (w) \, =\, \frac{D(z,w)}{(z - wq)(z-wq^{-1})} \quad |z|>|w|,
\end{equation}
where
\bray
\label{2dezao}
D(z,w) \, &=&\, \left( \frac{z}{w} \right)^{\sqrt{2} P} 
\exp \left\{ -\sqrt{2}\sum_{n<0}\frac{\alpha_n}{[n]}((zq^{-\frac{1}{2}})^{-n} 
\, -\, (wq^{\frac{1}{2}})^{-n}) \right\} \times \nonumber \\
&\times &
\exp \left\{ -\sqrt{2}\sum_{n>0}\frac{\alpha_n}{[n]}((zq^{\frac{1}{2}})^{-n} 
\, -\, (wq^{-\frac{1}{2}})^{-n}) \right\} \quad .
\eray          
The non-local structure of the poles in (\ref{2baker}) requires an expansion of
the numerator in q-Taylor series (see appendix). In order to get a 
symmetric result we  split the
numerator in two equal terms and expand each one separately around the diferent
poles. The OPE may be written as
\begin{equation}
E^{+} (z) E^{-} (w) \, \sim \, \frac{1}{w(q\, -\, q^{-1})} \left\{
\frac{\Psi (wq^{\frac{1}{2}})}{z\, -\, wq} \, -\, 
\frac{\Phi (wq^{-\frac{1}{2}})}{z\, -\, wq^{-1}} \right \} \,  ,\quad
|z|>|w| \, 
\label{ee}
\end{equation}
where 
\begin{equation}
\label{2psi}
\Psi (z) \, =\, q^{\sqrt{2} P} 
\exp \left\{ \sqrt{2} (q-q^{-1}) \sum_{n>0} \alpha_n z^{-n} \right \} \quad ,
\end{equation}
and
\begin{equation}
\label{2phi}
\Phi (z) \, =\, q^{-\sqrt{2} P} 
\exp \left\{ -\sqrt{2} (q-q^{-1}) \sum_{n<0} \alpha_n z^{-n} \right \}  \quad ,
\end{equation}
 where the symbol $\sim $ stand for  equality up to regular terms.
 The currents $E^{\pm} (z)$ are generating functions
of Drinfeld's generators for the q-deformed algebra 
${\cal{U}}_q ({\widehat{su}}(2))$, 
$ E^{\pm} (z) \, =\, \sum_{n \in \Bbb Z} E^{\pm}_n z^{-n-1}$
and the combination \\
$(\Psi (z) \, -\, \Phi (z))/(\sqrt{2} z(q\, -\, q^{-1}))$
correspond to  the q-deformed analogue of the undeformed Cartan sub algebra 
current $H(z)$ of ${\widehat{su}}(2)$. 
It should be noticed from (\ref{ee}) that the deformation naturally splits the
positive and negative oscillator modes according to positive and negative
powers of $q$.
 The full
 ${\cal{U}}_q ({\widehat{su}}(2))$ is obtained from (\ref{ee}) together with  
 \bray
\Psi (z) \Phi (w) \, &=&\, 
\frac{(z\, -\, wq^3 )(z\, -\, wq^{-3})}{(z\, -\, wq)(z\, -\, wq^{-1})} 
\Phi (w) \Psi (z) \quad , \nonumber \\
\Psi (z) E^{\pm} (w) \, &=&\, q^{\pm 2} 
\frac{(z\, -\, wq^{\mp \frac{5}{2}})}{(z\, -\, wq^{\pm \frac{3}{2}})}
E^{\pm}(w) \Psi (z) \quad , \nonumber \\
E^{\pm} (z) \Phi (w) \, &=&\,  q^{\pm 2} 
\frac{(z\, -\, wq^{\mp \frac{5}{2}})}{(z\, -\, wq^{\pm \frac{3}{2}})}
\Phi (w) E^{\pm} (z) \quad , \nonumber \\
E^{\pm} (z) E^{\pm} (w) \, &=&\, 
\frac{(zq^{\pm 2} \, -\, w)}{(z\, -\, wq^{\pm 2})} E^{\pm} (w) E^{\pm} (z)
\quad  |z|>|w|. 
\eray                                  

The regular part of the OPE in eq. (\ref{ee}) is important to define the
Energy Momentum tensor in the Sugawara form. Here we introduce the notation 
${\,}^{\times}_{\times} E^{+}(w)E^{-}(w)^{\times}_{\times}$ 
for the remaining regular terms after taking equal point limit, i.e.
\begin{equation}
 {\,}^{\times}_{\times}E^{+}(w)E^{-}(w)^{\times}_{\times} \, =\, 
\lim_{z \rightarrow w} \left\{ E^{+}(z)E^{-}(w) \, -\, 
 \frac{1}{w(q\, -\, q^{-1})} \left(
\frac{\Psi (wq^{\frac{1}{2}})}{z\, -\, wq} \, -\, 
\frac{\Phi (wq^{-\frac{1}{2}})}{z\, -\, wq^{-1}}  \right) \right\}
\end{equation}
  Due to the nonlocal character of
the q-Taylor expansion, it consist of an infinite number 
terms which vanish when $q \rightarrow 1$.
The first order term  is given by
\[
{\,}^{\times}_{\times}E^+ (w) E^- (w)^{\times}_{\times} \, =\,   
\frac{1}{2[2]w^2  (q\, -\, q^{-1})^2 } \left \{
[q^{-3} (A^{1}(wq^2 )  - 1)  + (q-q^{-1})]\Psi (wq^{\frac{1}{2}})
\right.  \,  
\]
\begin{equation}
+\, \left.
 \Phi (wq^{-\frac{1}{2}})[q^3 ( B^{1}(wq^{-2})-1)-(q-q^{-1})]
 \right \} \, + {\cal{O}}(q - q^{-1} )  \quad .
\end{equation}   
where the composite fields $A^{\alpha}$ and $B^{\alpha}$ are defined as
\begin{equation}
\label{azao}
A^{\alpha} (w) \, =\, \Phi^{-1} (wq^{-\frac{\alpha}{2}}) 
\Psi (wq^{\frac{\alpha}{2}}) \quad ,
\end{equation}
\begin{equation}
\label{bezao}
B^{\alpha} (w) \, =\, \Phi (wq^{-\frac{\alpha}{2}}) 
\Psi^{-1} (wq^{\frac{\alpha}{2}}) \quad .
\end{equation}
They are obtained by direct calculation, as explained in
appendix. We also prove that the higher power terms in $(q - q^{-1})$ 
appearing in the regular part of $E^+(z)E^-(w)$ only 
depend upon, either $D(wq^{2k+1},w)$ or $D(wq^{-2k-1},w)$ with their general
structure given by
\bray
\label{dezao}
D(wq^{2k+1},w) \, &=&\, :A^1 (wq^2) A^1 (wq^4) \ldots A^1 (wq^{2k}):
\Psi (wq^{\frac{1}{2}}) \quad , \nonumber \\
D(wq^{-2k-1},w) \, &=&\, \Phi (wq^{-\frac{1}{2}})
:B^1 (wq^{-2}) B^1 (wq^{-4}) \ldots B^1 (wq^{-2k}):\, .
\eray 
It should be noticed that all nonvanishing terms constituting 
 the regular part of $E^+ (z)E^- (z)$ are functionals of 
$\Psi $ and $\Phi $ and
henceforth, non linear functionals of the Cartan subalgebra current $H(z)$. This
fact suggests a generalization of the q-analogue of the quantum equivalence
theorem which establishes, in particular, 
 the equality of the energy momentum tensor constructed
out of a level 1 simply laced current algebra $\widehat{\frak{g}}$ and its 
Cartan subalgebra (see Goddard and Olive \cite{olive}).

\section{$N=2$ q-Superconformal Algebra}

The construction of the $ N=2$  superconformal algebra requires two 
supercharges, $G^{\pm}(z)$. The most singular term of their OPE is proporcional
to a triple pole, as we have seen in section 2 for the $N=1$ case.  
As we have argued before, the triple pole is expected to
be replaced by a product of simple poles symmetrically displaced in terms of
the deformation parameter $q$. Following the same line of reasoning, we define
the deformed supercharges as
\begin{equation}
G^{\pm} (z) \, =\, E^{\pm} (z) \psi (z) \quad ,
\end{equation}
where $E^{\pm}(z)$  denote the currents related to 
${\cal{U}}_q ({\widehat{su}}(2))$ step operators, and $\psi (z)$ is the
fermi field defined by expression (\ref{umfermiope}).  In particular, 
for level one representations, we  realize them as the vertex operators
(\ref{2vertex}). The OPE for the two supercharges yields
\bray
\label{gg}
G^{+}(z)G^{-}(w) \, &=&\,   \left\{ \frac{1}{(z\, -\, w)(z\, -\, wq)} \, +\, 
\frac{q^{\frac{1}{2}}L^{1}(wq^{\frac{1}{2}})}{(z\, -\, wq)} \,  \right.
\nonumber \\
&+&\, \left.
\frac{1}{2}
\frac{(q^{-3}(A^{1}(wq^2) \,-\, 1)\, +\, (q\,-\,q^{-1}))}
{[2] w(q\, -\, q^{-1})(z\, -\, w)}\, +\, {\cal O} (q-q^{-1}) 
\right\} \frac{\Psi (wq^{\frac{1}{2}})}{w(q\, -\, q^{-1})} \, 
\nonumber \\
&-&\, 
\frac{\Phi (wq^{-\frac{1}{2}})}{w(q\, -\, q^{-1})} \left\{ 
\frac{1}{(z\, -\, w)(z\, -\, wq^{-1})} \, +\, 
\frac{q^{-\frac{1}{2}}L^{-1}(wq^{-\frac{1}{2}})}{(z\,-\, wq^{-1})} \,  \right.
\nonumber \\
&+&\, \left.
\frac{1}{2}\frac{(q^{3}(B^{1}(wq^{-2})\, -\, 1)\, -\, (q\, -\, q^{-1}))}
{[2] w (q\, -\, q^{-1})(z\, -\, w)} \, +\, {\cal O}(q-q^{-1})\right\} \quad ,
\eray
where $L^{\alpha }$ is the fermionic Energy-Momentum tensor defined in
(\ref{umemfermi}), $A^{\alpha}$ and $B^{\alpha}$ are defined in 
(\ref{azao}, \ref{bezao}). Again, the shift in the arguments of currents in the
r.h.s. requires an additional index to label an infinite family of supercharge
generators, we therefore define  
\begin{equation}
\label{3supercarga}
G^{\pm}_{\alpha} (z) \, =\, E^{\pm} (zq^{\pm \frac{\alpha}{2}}) 
\psi (zq^{\mp \frac{\alpha}{2}}) \quad .
\end{equation}
The OPE relations for the supercharge generators is given by
\bray
\label{3superesuper2}
G^{+}_{\alpha} (z) G^{-}_{\beta} (w) \, &=&\,  
\frac{{\,}^{\times}_{\times} E^{+} (wq^{\frac{\beta}{2} +\alpha}) 
E^{-} (wq^{-\frac{\beta}{2}})^{\times}_{\times}}
{(zq^{-\frac{\alpha}{2}} \, -\, wq^{\frac{\beta}{2}})} \, 
\nonumber \\
&+&\, 
\frac{[1-\alpha -\beta]q^{\frac{1+\beta -\alpha}{2}} 
L^{1-\alpha -\beta} (wq^{\frac{1-\alpha}{2}}) \Psi (wq^{\frac{1-\beta}{2}})}
{zq^{\frac{\alpha}{2}} \, -\, wq^{-\frac{\beta}{2} +1}} \, 
\nonumber \\
&+&\, 
\frac{[1+\alpha +\beta]q^{\frac{1+\beta +\alpha}{2}} 
\Phi (wq^{\frac{-1-\beta}{2}}) L^{1+\alpha +\beta} (wq^{\frac{1+\alpha}{2}})} 
{zq^{\frac{\alpha}{2}} \, -\, wq^{-\frac{\beta}{2} -1}} \,  
\nonumber \\
&+&\,  
\frac{1}{wq^{-\frac{\beta}{2}}(q\, -\, q^{-1})
(zq^{-\frac{\alpha}{2}} \, -\, wq^{\frac{\beta}{2}})} \left\{
\frac{\Psi (wq^{\frac{1-\beta}{2}})}
{zq^{\frac{\alpha}{2}} \, -\, wq^{-\frac{\beta}{2} +1}} \,  \right.
\nonumber \\
&-&\, \left.
\frac{\Phi (wq^{\frac{-1-\beta}{2}})}
{zq^{\frac{\alpha}{2}} \, -\, wq^{-\frac{\beta}{2} -1}} \right\} \quad ,
\eray
where
\bray
\label{regpart}
{\,}^{\times}_{\times}E^+ (wq^{\frac{\beta}{2} +\alpha}) 
E^- (wq^{-\frac{\beta}{2}})^{\times}_{\times} \, &=&\,   
\frac{1}{2[2]w^2 q^{-\beta}  (q\, -\, q^{-1})^2 } \left \{
[q^{-3} (A^{1}(wq^{-\frac{\beta}{2}+2})  - 1)  \, \right.
\nonumber \\ 
&+&\, (q-q^{-1})] \Psi (wq^{\frac{-\beta +1}{2}})  \,  
\nonumber \\
&+&\, 
\Phi (wq^{\frac{-\beta -1}{2}})
[q^{3} (B^{1}(wq^{-\frac{\beta}{2}-2}) - 1) \, 
\nonumber \\
&-&\, \left. 
(q-q^{-1})] \right\} \, + \, {\cal{O}} (q-q^{-1}) \quad .
\eray

In the classical limit ($q\rightarrow 1$) these generators
$G^{\pm}_{\alpha}$ become the usual supercharges of $N=2$ superconformal
algebra,  $G^{\pm} (z) \, =\, E^{\pm} (z) \psi (z)$ obeying the OPE relation,
(see \cite{watterson}) 
\[
G^{+} (z) G^{-} (w) \, =\, \frac{1}{(z\, -\, w)^{3}} \, +\, 
\frac{\sqrt{2} H(w)}{(z\, -\, w)^{2}} \, +\, 
\frac{2T(w)\, +\, \sqrt{2} H'(w)}{(z\, -\, w)} \quad ,
\]
where  $T(w)=T_{bosonic} (w) +L_{fermionic} (w)$ is the total Energy-Momentum
tensor, as in $N=1$ case. 
Analysing expressions (\ref{3superesuper2}) and (\ref{regpart}), we
can formulate an expression for the bosonic Energy-Momentum tensor taking the
coefficient of simple poles in the OPE, bearing in mind the correct classical
limit to be \[
T(z)\, =\, \frac{1}{2} :H^{2} (z): \quad .
\]
We therefore define a family of operators labeled by two indices 
\begin{equation}
\label{3qenermom}
T^{\alpha ;\beta} (z) \, =\, \frac{1}{2[2]z^2 (q\, -\, q^{-1})^2} \left\{
A^{\alpha} (zq^{\beta}) \, +\, B^{\alpha} (zq^{\beta}) \, -\, 2 \right\}
\end{equation}
with $A^{\alpha}$ and $B^{\alpha}$ defined by expressions
(\ref{azao},\ref{bezao}). The algebra of the q-deformed Energy-Momentum tensor
(\ref{3qenermom}) does not close, however its components
$A^{\alpha}$ and $B^{\alpha}$ define a quadratic algebra, namely
\begin{eqnarray}
\label{3superalgebra}
A^{\alpha} (z) A^{\beta} (w) \, &=&\, 
f(zq^{\frac{\alpha}{2}};wq^{-\frac{\beta}{2}}) :A^{\alpha} (z) A^{\beta} (w):
\quad , \nonumber \\
B^{\alpha} (z) B^{\beta} (w) \, &=&\, 
f(zq^{\frac{\alpha}{2}};wq^{-\frac{\beta}{2}}) :B^{\alpha} (z) B^{\beta} (w):
\quad , \nonumber \\
A^{\alpha} (z) B^{\beta} (w) \, &=&\, 
f^{-1} (zq^{\frac{\alpha}{2}};wq^{-\frac{\beta}{2}}) 
:A^{\alpha} (z) B^{\beta} (w):
\quad , \nonumber \\
B^{\alpha} (z) A^{\beta} (w) \, &=&\, 
f^{-1} (zq^{\frac{\alpha}{2}};wq^{-\frac{\beta}{2}}) 
:B^{\alpha} (z) A^{\beta} (w):
\quad , \nonumber \\
A^{\alpha} (z) G^{\pm}_{\beta} (w) \, &=&\, 
q^{\pm 2} 
\frac{(zq^{\frac{\alpha}{2}}\, -\, wq^{\pm \left( \frac{\beta -5}{2} \right)})}
{(zq^{\frac{\alpha}{2}}\, -\, wq^{\pm \left( \frac{\beta +3}{2} \right)})}
:A^{\alpha} (z) G^{\pm}_{\beta} (w):
\quad , \nonumber \\
B^{\alpha} (z) G^{\pm}_{\beta} (w) \, &=&\, 
q^{\mp 2} 
\frac{(zq^{\frac{\alpha}{2}}\, -\, wq^{\pm \left( \frac{\beta +3}{2} \right)})}
{(zq^{\frac{\alpha}{2}}\, -\, wq^{\pm \left( \frac{\beta -5}{2} \right)})}
:B^{\alpha} (z) G^{\pm}_{\beta} (w):
\quad , \nonumber \\
L^{\alpha} (z) G^{\pm}_{\beta} (w) \, &=&\, 
\frac{1}{[\alpha] w(q\, -\, q^{-1})} \left\{
\frac{G^{\pm (\pm \beta -\alpha)} (wq^{\frac{\alpha}{2}})}
{z\, -\, wq^{\mp \frac{\beta}{2} + \frac{\alpha}{2}}} \, + \right. \nonumber \\
&-&\,\left.
\frac{G^{\pm (\pm \beta +\alpha)} (wq^{\frac{-\alpha}{2}})}
{z\, -\, wq^{\mp \frac{\beta}{2} - \frac{\alpha}{2}}} \ \right\} \quad ,
\end{eqnarray}
where the analytic structure is given by
\[
f(z;w) \, =\, 
\frac{(z\, -\, wq)(z\, -\, wq^{-1})}{(z\, -\, wq^{3})(z\, -\, wq^{-3})} \quad ,
\]
and $L^{\alpha} (z)$ is the fermionic Energy-Momentum tensor as defined by
expression (\ref{umemfermi}).

The higher order terms in the q-Taylor expansion of the OPE (see appendix) lead
us to a more general class of q-deformed Energy-Momentum tensors in Sugawara
form. In fact, the multi index family of operators given by
\begin{eqnarray}
\label{3sugawara}
T^{\alpha_1 ,\beta_1 ; \ldots ;\alpha_n ,\beta_n} (z) \, &=&\,
\frac{1}{2 n^2 [2] z^2 (q\, -\, q^{-1})^2} \left\{
:A^{\alpha_1} (zq^{\beta_1}) \cdots A^{\alpha_n} (zq^{\beta_n}): \, +
\right. \nonumber \\
&+&\, \left.
:B^{\alpha_1} (zq^{\beta_1}) \cdots B^{\alpha_n} (zq^{\beta_n}): \, -\, 2 
\right\} \quad ,
\end{eqnarray}
with dots indicating that all $\Phi$ are placed on the left of $\Psi$. As
before, the algebra of the multi index Energy-Momentum tensor
(\ref{3sugawara}) does not close however their constituents.
\begin{equation}
\label{A}
A^{\alpha_1 ,\beta_1 ;\ldots ;\alpha_n ,\beta_n} (z) \, =\,
:A^{\alpha_1} (zq^{\beta_1}) \cdots A^{\alpha_n} (zq^{\beta_n}): 
\end{equation}
and
\begin{equation}
\label{B}
B^{\alpha_1 ,\beta_1 ;\ldots ;\alpha_n ,\beta_n} (z) \, =\,
:B^{\alpha_1} (zq^{\beta_1}) \cdots B^{\alpha_n} (zq^{\beta_n}): 
\end{equation}
close two isomorphic OPE algebras of the type:
\[
A^{\alpha_1 ,\beta_1 ;\ldots ;\alpha_n ,\beta_n} (z)
A^{\gamma_1 ,\delta_1 ;\ldots ;\gamma_m ,\delta_m} (w) \, = 
\]
\begin{equation}
\label{3azaoalgebra}
=\, \prod_{i=1}^{n} \prod_{j=1}^{m} f(zq^{\beta_i +\frac{\alpha_i}{2}} ; 
wq^{\delta_j -\frac{\gamma_j}{2}})
:A^{\alpha_1 ,\beta_1 ;\ldots ;\alpha_n ,\beta_n} (z)
A^{\gamma_1 ,\delta_1 ;\ldots ;\gamma_m ,\delta_m} (w): \quad .
\end{equation}

Notice that the OPE  $G^+ (z)G^- (w)$
naturally produces terms of the type (\ref{A}) and (\ref{B}), because  
$D(wq^{2k+1},w) = A^{1,2;1,4;...;1,2k} (w) \Psi (wq^{\frac{1}{2}})$  and \\
$D(wq^{-2k-1},w) = \Phi (wq^{-\frac{1}{2}}) B^{1,-2;1,-4;...;1,-2k}(w)$
 given in (\ref{dezao}).  These are particular cases when $\alpha_i = 1$,
$\beta_i =2i$ in the case of $A$, and $\beta_i =-2i$ in the case of $B$, 
leading to a remarkable
cancellation of poles in (\ref{3azaoalgebra}) in which the $2mn$ poles result
in $2m$ only.

\section{Sugawara Construction and ${\cal{U}}_{q} ({\widehat{su}}(N+1))$ Vertex
Operators}

The vertex operator construction for an affine simply laced Lie algebra \\
${\cal{U}}_q ({\widehat{su}}(N+1))$, requires  N  independent copies of 
q-oscillators satifying 
\begin{equation}
\label{4osc}
[\alpha^{i}_{n} ,\alpha^{j}_{m}] \, =\, \frac{[K_{ij} n][n]}{2n} \delta_{m+n,0}
\; , i=1\ldots N \quad ,
\end{equation}
where $K_{ij}$ is the Cartan matrix of the underlying Lie algebra $su(N+1)$, and
 zero modes $P^{i}$ and $Q^{i}$ satisfying usual commutation
relations
\[
[Q^{i} ,P^{j} ] \, =\, i \delta^{ij} \quad .
\]
These commutation relations become, in the limit $q\rightarrow 1$,
 the usual ones for the modes of the
Cartan sub-algebra of the undeformed ${\widehat{su}}(N+1)$ affine Lie algebra
in the Chevalley basis, i.e.
\[
[\alpha^{i}_{n} ,\alpha^{j}_{m} ] \, =\, K_{ij}n \, \delta_{m+n,0} \quad .
\]

Let $\beta_i \; i=1\ldots N$ be the simple roots of $g$ algebra. To each
simple root we define vertex operators as indicated in ref. \cite{fjing}
\bray
\label{4simpleroots}
E^{\pm}_{i} (z)\, &=&\, \exp \left\{ \mp \sqrt{2} \sum_{n<0}
\frac{\alpha^{i}_{n}}{[n]} (zq^{\mp \frac{1}{2}})^{-n} \right\}\, \times
\nonumber \\
&\, &\times \, \exp \left\{ \mp \sqrt{2} \sum_{n>0} \frac{\alpha^{i}_{n}}{[n]}
(zq^{\pm \frac{1}{2}})^{-n} \right\}e^{\pm i \sqrt{2} Q^{i}} 
z^{\pm \sqrt{2} P^{i}} \quad , \nonumber \\
\Psi_{i} (z) \, &=&\, q^{\sqrt{2} P^{i}} \exp \left\{ \sqrt{2} 
(q\, -\, q^{-1}) \sum_{n>0} \alpha^{i}_{n} z^{-n} \right\} \quad , \nonumber \\
\Phi_{i} (z) \, &=&\, q^{-\sqrt{2} P^{i}} \exp \left\{ -\sqrt{2} 
(q\, -\, q^{-1}) \sum_{n<0} \alpha^{i}_{n} z^{-n} \right\} \quad , 
\eray
satisfying  the OPE relations
\bray
\label{4operelations}
E^{+}_{i} (z) E^{-}_{i} (w) \, &\sim & \frac{1}{w(q\, -\, q^{-1})}
\left\{ \frac{\Psi_{i} (wq^{\frac{1}{2}})}{z\, -\, wq^{1}} \, -\, 
\frac{\Phi_{i} (wq^{-\frac{1}{2}})}{z\, -\, wq^{-1}} \right\} \quad , 
\nonumber \\
E^{-}_{i} (z) E^{+}_{i} (w) \, &\sim & \frac{1}{w(q\, -\, q^{-1})}
\left\{ \frac{\Phi_{i} (wq^{\frac{1}{2}})}{z\, -\, wq^{1}} \, -\, 
\frac{\Psi_{i} (wq^{-\frac{1}{2}})}{z\, -\, wq^{-1}} \right\} \quad , 
\nonumber \\
E^{+}_{i} (z) E^{-}_{j} (w) \, &=& \frac{(z\, -\, w)}{z} 
:E^{+}_{i} (z) E^{-}_{j} (w): \; , i=j\pm 1 \quad , 
\nonumber \\
\Psi_{i} (z) E^{\pm}_{i} (w) \, &=&\, q^{\pm 2} 
\frac{(z\, -\, wq^{\mp \frac{5}{2}})}{(z\, -\, wq^{\pm \frac{3}{2}})}
E^{\pm}_{i} (w) \Psi_{i} (z)\quad , 
\nonumber \\
\Psi_{i} (z) E^{\pm}_{j} (w) \, &=&\,  
\frac{(z\, -\, wq^{\pm \frac{1}{2}})}{(z\, -\, wq^{\mp \frac{3}{2}})}
E^{\pm}_{j} (w) \Psi_{i} (z)\; i=j\pm 1 \quad , 
\nonumber \\
E^{\pm}_{i} (z) \Phi_{i} (w) \, &=&\, q^{\pm 2} 
\frac{(z\, -\, wq^{\mp \frac{5}{2}})}{(z\, -\, wq^{\pm \frac{3}{2}})}
\Phi_{i} (w) E^{\pm}_{i} (z)  \quad , 
\nonumber \\
E^{\pm}_{i} (z) \Phi_{j} (w) \, &=&\,  
\frac{(z\, -\, wq^{\pm \frac{1}{2}})}{(z\, -\, wq^{\mp \frac{3}{2}})}
\Phi_{j} (w) E^{\pm}_{i} (z)  \; i=j\pm 1 \quad , 
\nonumber \\
\Psi_{i} (z) \Phi_{i} (w) \, &=& 
\frac{(z\, -\, wq^3 )(z\, -\, wq^{-3})}{(z\, -\, wq)(z\, -\, wq^{-1})}
\Phi_{i} (w) \Psi_{i} (z) \quad , 
\nonumber \\
\Psi_{i} (z) \Phi_{j} (w) \, &=& 
\frac{(z\, -\, w)^{2}}{(z\, -\, wq^{2})(z\, -\, wq^{-2})}
\Phi_{j} (w) \Psi_{i} (z) \; i=j\pm 1 \quad .
\eray

For a positive root  $\beta \, =\, \beta_{i}
+\beta_{i+1}+\cdots +\beta_{i+s} \; 0\leq s\leq N-i$, we consider the 
 normal ordered  products \cite{fjing}
\bray
\label{4nonsimpleroots}
E^{+}_{\beta} (z) \, &=&\, :E^{+}_{i} (zq^{i}) E^{+}_{i+1} (zq^{i+1}) \ldots
E^{+}_{i+s} (zq^{i+s}) : \quad , 
\nonumber \\
E^{-}_{-\beta} (z) \, &=&\, :E^{-}_{i+s} (zq^{i+s}) E^{-}_{i+s-1} (zq^{i+s-1}) 
\ldots E^{-}_{i} (zq^{i}) : \quad , 
\nonumber \\
\Psi_{\beta} (z) \, &=&\, \Psi_{i} (zq^{i}) \Psi_{i+1} (zq^{i+1}) \ldots
\Psi_{i+s} (zq^{i+s}) \quad ,
\nonumber \\
\Phi_{\beta} (z) \, &=&\, \Phi_{i} (zq^{i}) \Phi_{i+1} (zq^{i+1}) \ldots
\Phi_{i+s} (zq^{i+s}) \quad ,
\eray
The OPE relations for these
operators are obtained using rules (\ref{4operelations}). In particular the 
OPE of  $E^{+} (z) E^{-} (w)$ we find
\[
\zeta (q,z)E^{+}_{\beta} (z) E^{-}_{-\beta} (w) \, =\, 
\frac{1}{w(q\, -\, q^{-1})}
\left\{ \frac{\Psi_{\beta} (wq^{\frac{1}{2}})}{z\, -\, wq^{1}} \, -\, 
\frac{\Phi_{\beta} (wq^{-\frac{1}{2}})}{z\, -\, wq^{-1}} \right\} \,  
\]
\[
+\, \frac{1}{2[2]w^{2}q^{2}(q\, -\, q^{-1})^{2}} \left\{ [q^{-3} 
(A^{1}_{\beta} (wq^{2}) \, -\, 1) \, + (q\, -\, q^{-1})] 
\Psi_{\beta} (wq^{\frac{1}{2}}) \right.
\]
\begin{equation}
\label{4nonsimplevertex}
+\, \left. \Phi_{\beta} (wq^{-\frac{1}{2}})
[ q^{3} (B^{1}_{\beta} (wq^{-2}) \, -\, 1) \, -\, (q\, -\, q^{-1})] \right\}  
\, + {\cal{O}} (q\, -\, q^{-1}) \quad , 
\end{equation}
where
\[
\zeta (q,z) \, =\, z^{2(s-1)} q^{i(2s-1) +s(s+1)} \quad ,
\]
is a normalizer term introduced only to cancel factors apearing in OPEs of
vertex associated to neighbouring roots. The fields $A^{k}_{\beta}$ and
$B^{k}_{\beta}$ are defined as
\bray
A^{k}_{\beta} (z) \, &=&\,  
 :A^{k}_{i} (zq^{i}) A^{k}_{i+1} (zq^{i+1}) \ldots A^{k}_{i+s} (zq^{i+s}):
\quad , 
\nonumber \\
B^{k}_{\beta} (z) \, &=&\, 
 :B^{k}_{i} (zq^{i}) B^{k}_{i+1} (zq^{i+1}) \ldots B^{k}_{i+s} (zq^{i+s}):
\quad .
\eray

The Energy-Momentum tensor associated to ${\cal U}_{q} ({\widehat{su}}(N+1))$
is therefore proposed to be
\bray                                                              
\label{4qsugawara}
T^{k,l} (z) \, &=&\, \frac{1}{2(2\, +\, N)[2]z^2 (q\, -\, q^{-1})^{2}} \left\{
\sum_{i=1}^{N} [A^{k}_{i} (zq^{l}) \, +\, 
B^{k}_{i} (zq^{l}) \, -\, 2] \, + \right. 
\nonumber \\
&+&\, \left.
2\sum_{\beta >0} [A^{k}_{\beta} (zq^{l}) \, +\, 
B^{k}_{\beta} (zq^{l}) \, -\, 2]
\right\} \quad ,
\eray
providing a closed quadratic algebra for its constituents $A^{k}_{\beta}$ and 
$B^{k}_{\beta}$ which can be obtained using the primitive OPE relations in
(\ref{4operelations}). In the limit $q\rightarrow 1$, we recover the usual
Sugawara construction, namely \cite{olive} 
\begin{equation}
\label{4classenermom}
T(z) \, =\, \frac{1}{2(1+h)} \left\{ \sum_{i=1}^{N-1} {\,}^{\times}_{\times} 
H_{i}^{2} (z)^{\times}_{\times} \, +\, \sum_{\beta >0} {\,}^{\times}_{\times}
E^{\beta} (z) E^{-\beta} (z) \, +\, 
E^{-\beta} (z) E^{\beta} (z)^{\times}_{\times} \right\}\quad ,
\end{equation}
where $h$ is the dual Coxeter number (for the algebra $su(N+1)$ we have 
$h=N+1$),
and the crosses indicate normal ordering for the modes of $su(N+1)$ 
currents. For level 1 simply laced algebras, the contribution of terms
quadratic in step operators in the Energy-Momentum tensor is known to be 
proportional to terms dependent on the Cartan sub-algebra, see \cite{olive}
\[
{\,}^{\times}_{\times} E^{\beta} (z) E^{-\beta} (z) \, +\, 
E^{-\beta} (z) E^{\beta} (z)^{\times}_{\times} \, =\, 
{\,}^{\times}_{\times} (\beta . H(z))^{2}{\,}^{\times}_{\times} \quad .
\]

\section{Conclusion and Outlook}

A q-deformed version of the $N=2$ superconformal algebra was proposed in terms
of level 1 representations of ${\cal U}_{q} ({\widehat{su}}(2))$ Kac-Moody
algebra and a single real fermi field.  This construction hints the form for
the q-analogue for
 the Sugawara Energy-Momentum tensor possessing an exponential
dependence on the Cartan subalgebra generators.  An interesting point to be
further investigated concerns the construction of an action from which  such
Energy-Momentum tensor could be obtained using canonical methods following the
line of ref \cite{devega}.  

The study of the representations of such algebraic structure also deserves to
be further developed.  In particular, identities involving conformal embeddings
\cite{ago} as well as coset constructions \cite{gko} may be obtained from
representation theory.

{\large {\bf Aknowledgements}:}
 We thank Profs. G. M. Sotkov and A. H. Zimerman for many
helpful discussions. One of us (JFG) thanks  the International Centre for 
Theoretical Physics (Trieste) for hospitality
and support where part of this work was done.

\appendix 
\section{Expansions in q-Taylor Series}

In this appendix we review some properties of the q-Taylor expansion, which is
obtained from usual Taylor series by adding and subtracting terms. This
ensures the same content of the series in classical analysis, with the
additional advantage of providing a  non local expansion. An analytic
function $f(z)$ may be written near $z=\omega$ either  

\begin{equation}
\label{apqtaylor1}
f(z) \, =\, \sum_{k=0}^{\infty}  \frac{\partial^{2k}_{q} f(wq)}{[2k]!} 
(z\, -\, w)^{2k}_{q} \, +\,  \sum_{k=0}^{\infty}  
\frac{\partial^{2k+1}_{q} f(w)}{[2k+1]!} (z\, -\, wq)^{2k+1}_{q} \quad ,
\end{equation}
or  
\begin{equation}
\label{apqtaylor2}
f(z) \, =\, \sum_{k=0}^{\infty}  \frac{\partial^{2k}_{q} f(wq^{-1})}{[2k]!} 
(z\, -\, w)^{2k}_{q} \, +\,  \sum_{k=0}^{\infty}  
\frac{\partial^{2k+1}_{q} f(w)}{[2k+1]!} (z\, -\, wq^{-1})^{2k+1}_{q} \quad ,
\end{equation}
where $[n]! =[n][n-1]!$, the symbol $\partial_q$ denotes the q-derivative 
defined as
\begin{equation}
\label{apqderivada}
\partial_{q} f(z)\, =\, \frac{f(zq) \, -\, f(zq^{-1})}{z(q\, -\, q^{-1})} 
\quad , 
\end{equation}
and 
\begin{equation}
\label{apqbinomio}
(z\, -\, w)^{n}_{q} \, =\, \prod_{k=1}^{n} (z\, -\, wq^{n-2k+1}) \, =\, 
\sum_{k=0}^{n} \frac{[n]!}{[k]! [n-k]!} z^{k} (-w)^{n-k} \quad ,
\end{equation}
is the q-binomial. These expansions display a non local character and 
become useful
in  finding regular parts of OPE in q-deformed meromorphic field 
theories like those presented in this paper.

From equation (\ref{apqderivada}) we propose a general closed expression
 for the n-th q-derivative 
\begin{equation}
\label{apnqderivada}
\partial_{q}^{n} f(z) \, =\, \frac{1}{z^{n} (q\, -\, q^{-1})^{n}}
\sum_{\sigma_{0} =\pm 1} \sum_{\sigma_{1} =\pm 1} \cdots 
\sum_{\sigma_{n-1} =\pm 1} \left( \prod_{i=0}^{n-1} \sigma_{i} q^{-i\sigma_{i}}
\right) f\left( zq^{\sum_{i=0}^{n-1} \sigma_{i}} \right) \quad ,
\end{equation}
which can be easily proved by induction.

A crucial observation from (\ref{apnqderivada}) is that the $n^{th}$ order 
q-derivative is proportional to the original function with the argument 
shifted by some power of the deformation parameter.  The n sums over 
$\sigma _i = \pm 1$, $i=0, \cdots n-1 $ leads to a shifted argument 
proportional to $q^{n-2\alpha}$.  Henceforth, the parity of the power of $q$ 
is always the same
as the order of the derivative.  On the other hand, odd derivative terms in the
q-Taylor expansion (\ref{apqtaylor1}, \ref{apqtaylor2}) are evaluated in $z=w$
while the even ones in $z=wq$ or $z=wq^{-1}$.  In either case, the expansion of
$f(z)$ in (\ref{apqtaylor1}, \ref{apqtaylor2}) consist of linear combination of
terms like $f(wq^{2k+1}) \quad ; \, k\in \Bbb Z$.  

We now apply this result to calculate the regular part of the OPE 
$E^+(z)E^-(w)$ in eq. (\ref{2baker}).  The above argument implies that the 
r.h.s.  of (\ref{2baker})
is a linear combination of $D(wq^{2k+1},w)$.  The peculiar form of $D(z,w)$
yields a dependence entirely in terms upon the Cartan subalgebra current $H(z)$.
     To be more specific, the following identities can be obtained by direct
calculation. 
\bray
\label{appositiva}
D(wq^{2k+1},w) \,&=& \, \left( \prod_{i=0}^{k-1} 
\Phi^{-1} (wq^{-\frac{1}{2} +2(k-i)}) \right) \left( \prod_{j=0}^{k} 
\Psi (wq^{\frac{1}{2} +2(k-j)}) \right) \, = 
\nonumber \\
&=&\, 
:A^{1} (wq^{2}) A^{1} (wq^{4}) \ldots A^{1} (wq^{2k}): \Psi (wq^{\frac{1}{2}})
\quad .
\eray
and 
\bray
\label{apnegativa}
D(wq^{-(2k+1)},w) \,&=& \, \left( \prod_{i=0}^{k} 
\Phi (wq^{-\frac{1}{2} -2i}) \right) \left( \prod_{j=0}^{k-1} 
\Psi^{-1} (wq^{\frac{1}{2} -2j}) \right) \, = 
\nonumber \\
&=&\, 
\Phi (wq^{-\frac{1}{2}}) 
:B^{1} (wq^{-2}) B^{1} (wq^{-4}) \ldots B^{1} (wq^{-2k}):
\, .
\eray

We should also point out that due to the nonlocality, the expansion of the
regular part of $E^+(z)E^-(w)$ does not truncate when equal point limit is
taken.  It is, in fact composed of infinite number of terms which are
classified in powers of $(q-q^{-1})$.


\begin{thebibliography}{99}

\bibitem{afgz}
H. Aratyn, L.A. Ferreira, J.F. Gomes and A.H. Zimerman: Phys. lett. B254,
(1991) 372.

\bibitem{ago}
R. Arcuri, J.F. Gomes and D. Olive: Nucl. Phys. B285 (1987) 327


\bibitem{boug2}
A. H. Bougourzi and L. Vinet: Int. J. Mod. Phys. A10 (1995)  923.

\bibitem{cgk}
M. Chaichian, J.F. Gomes and P. Kulish: Phys. Lett. B311 (1993) 93.

\bibitem{chai}
M. Chaichian, P. Presnajder: Phys. Lett. B277 (1992) 109-118.

\bibitem{cgr}
M. Chaichian, J.F. Gomes and R. Gonzalez-Felipe: Phys. Lett. B341 (1994)
147-152.

\bibitem{devega}
H. De Vega and N. Sanchez: Phys. Lett. B216 (1989) 97-102.

\bibitem{fgsz}
L.A. Ferreira, J.F. Gomes, A. Schwimmer and A.H. Zimerman: Phys. Lett. B274,
(1992),65.        

\bibitem{fr}
E. Frenkel and N. Reshetikhin: ``Quantum Affine Algebras and 
Deformations of the Virasoro and W-algenbras'' preprint q-alg/9505025.

\bibitem{fjing}
I. B. Frenkel and N. H. Jing: Proc. Nat'l. Acad. Sci. (USA) 85(1988) 9373.

\bibitem{gko}
P. Goddard, A. Kent and D. Olive: Comm. in Math. Phys. 103 (1986) 105-119.

\bibitem{olive}
P. Goddard and D. Olive: Int. J. Mod. Phys. A1 (1986) 303-414.

\bibitem{goddard}
P. Goddard and A. Schwimmer: Phys. Lett. B214 (1986) 209-214.

\bibitem{singh}
C.H. Oh and K. Singh: ``Realizations of the q-Heisenberg and 
q-Virasoro algebras'' preprint hep-th/9408001.

\bibitem{sato1}
H. Sato: Nucl. Phys. B393 (1993) 442-458.

\bibitem{sato2}
H. Sato: preprint ``OPE Formulae for Deformed Super Virasoro Algebras ''
 NBI-HE-95-34 (1995).

\bibitem{watterson}
G. Waterson: Phys. Lett B171 (1986) 77.
               
\end{thebibliography}
\end{document}